# NOVEL SPANNING-TREE MATRIX APPROACH TO MODEL AND OPTIMIZE LARGE-SCALE, TREE-SHAPED WATER DISTRIBUTION NETWORKS

K. H. M. R. N. Senavirathna[1,*] and C. K. Walgampaya[1]

[1]Department of Engineering Mathematics, Faculty of Engineering, University of Peradeniya, Peradeniya, Sri Lanka

*Corresponding Author's email address: rajithas@eng.pdn.ac.lk

## ABSTRACT

There exist many criteria for the optimal design of water distribution networks. One of the most common criteria is to design the optimal cost water distribution network while satisfying the hydraulic design constraints. This study was carried out to propose a novel computational method named Spanning-Tree Matrix Approach that can model large-scale tree-shaped water distribution networks. A case study was tested to demonstrate the use of the Spanning-Tree Matrix Approach model coupled with the Honey-Bee Mating Optimization algorithm to find the combination of pipe diameters that minimizes the cost of the network. The results show that the Spanning-Tree Matrix Approach is successful in modeling a tree-shaped water distribution network of any size. Moreover, proposed Spanning-Tree Matrix Approach has the flexibility to be adapted to any desirable governing equation or design criteria being imposed, and the element of simplicity to output desired constraint evaluations into a modern stochastic optimization algorithm (i.e., Genetic Algorithm, Simulated Annealing, Ant-Colony Optimization, Honey-Bee Mating Optimization, etc.) for the network optimization purpose.

# INTRODUCTION

The recent past has witnessed many developments in computational methods in water resources management, with particular interest in water distribution networks (WDNs) modeling and optimization. The influence of evolving mathematics upon these developments of WDNs is inevitable. Rapid urbanization, scarcity of natural resources, and upscale-costs involved are identified to be the major reasons for this continuous research interest on WDNs. Among many criteria for the design of WDNs, one of the widely-used is, finding the combination of pipe diameters that satisfies the minimum-nodal-hydraulic-head and maximum-pipe-head-loss requirements while minimizing the cost involved. Hence, this is regarded as a combinatorial optimization problem. Due to the environmental conditions, the non-linear nature of the information associated with the hydraulic constrains, and the pipe materials and the type of architecture being adopted, the modeling of a given WDN has, most of the time, design constraints specific to its own. When existing hydraulic modeling packages do not facilitate the ability to model and output adaptable constraint evaluations to be used by optimization algorithms, the development of problem specific new computational methods by engineers and analysts becomes essential.

# LITERATURE REVIEW

A growing body of research shows that researchers have been using numerous mathematical methods to model and optimize water systems (Alporevits and Shamir, 1997; Lin et al., 1997; Cunha and Sousa, 1999; Vamvakeridou-Lyroudia et al., 2007; Todini and Rossman, 2013; Awe et al., 2020).

Huddleston et al. 2004 have brought about a clear, yet descriptive analysis of WDNs using one of the commonly used spreadsheets, *Microsoft Excel*. Extending the analysis and engineering design, there has been many computer programs to model WDNs, and they are being used by practitioners and researchers at present. As far as modeling of WDN is concerned, it is very noticeable that the use of a commercially available hydraulic modeling software package is still very common practice among the researchers, case in point, the use of EPANET by Mohan and Babu 2010, Suribabu 2012, Awe et al. 2020. It is well-known that commercially available hydraulic modeling software packages of this sort provide users considerable convenience when it comes to modeling a WDN. However, constraint evaluations in those packages may give off undesired, or rather impractical, results upon simulation, given the type of constraint evaluation being considered for optimization purpose. Therefore, it has been an apparent need to keep researching for novel computational methods that are highly adaptable to specific governing equations, and can output desired constraint evaluations of the WDN problem being concerned, given a reasonable complexity of the problem.

Although the incorporation of graph theory to tackle WDN has been rarely noticeable in the distant past, the case is not the same in the present context. The work done by Cong and Zhao 2015 has focused on applying a minimum spanning-tree algorithm, the Kruskal's algorithm, to find the minimum spanning-tree for a WDN that has fixed nodal positioning. Kruskal's algorithm, however, finds the minimum-weight spanning-tree by choosing edges with minimum weights while avoiding potential loops. Their study has assigned constant weights to the edges, or rather pipes, not taking into account the possibility of using commercially available pipes of different sizes. On the contrary to their study, it also needs to be noted that designing and optimizing WDNs should be done taking hydraulic constraints into consideration, as final goal is, presumably, the satisfaction of the consumer demand constraints.

Corte and Sörensen 2013, in their critical review, have shown that not only modeling, but also optimization of WDNs is much needed to reach the goals of delivering reliable water supply at optimum cost. In literature, there can be seen three main categories for WDN optimization algorithms, namely, linear-integer-programing, deterministic-heuristic algorithms and stochastic-heuristic algorithms.

In the study done by Samani and Mottaghi 2006, WDN modeling and optimization has been carried out with the use of coupled hydraulic and network optimization, where linear-integer-programming has been used as the optimization technique. In this method, they have turned the non-linear constraints of the optimization problem into linear form. Although this method has many advantages, such as simplicity and convergence, it also has disadvantages including the error when transforming the non-linear problem into linear one.

Deterministic-heuristic algorithms provide methods incorporating gradient and implicit information associated with the WDN to reach an optimal solution. The main advantage of algorithms of this kind is that, it takes a very smaller number of iterations to reach the optimal solution compared to stochastic-heuristic algorithms. Moreover, deterministic-heuristic algorithms cannot guarantee if the solution so obtained is the global optimum solution, as it can also be a local optimum. The studies done by Lin et al. 1997, Hsu and Cheng 2002, Mohan and Babu 2009, Suribabu 2012, and Awe et al. 2020, provide criteria for employing deterministic-heuristic algorithms to reach near optimal solutions for WDNs.

On the other hand, stochastic-heuristic, or sometimes called modern-stochastic-heuristic, methods neither take gradient nor implicit information associated with the WDNs. These kinds of algorithms evaluate the objective function at randomly taken different regions of the solution space to seek the feasible solution that minimizes the cost function. The main advantage of this method, as described earlier, is that it investigates the whole solution space to probe the global optimum solution. On the contrary to deterministic-heuristic methods, stochastic-heuristic methods perform extremely large number of objective function evaluations to reach this global optimum. Some of the notable studies that utilize these optimization algorithms for WDNs include Simulated Annealing (SA) Approach done by Cunha and Sousa 1999, Genetic Algorithm (GA) by Dijk et al. 2008, Differential Evolution (DE) by Moosavian and Lence 2019, and Honey-Bee Mating Optimization (HBMO) by Mohan and Babu 2010, and Senavirathna et al. 2022. More recently, El-Ghandour and Elbeltagi 2018, have carried out a remarkable study employing five stochastic-heuristic, evolutionary algorithms – GA, Particle Swarm Optimization (PSO) Algorithm, Ant Colony Optimzation (ACO) Algorithm, Memetic Algorithm (MA), and Shuffled Frog Leaping Algorithm (SFLA) - to obtain the optimal design of two benchmark WDNs.

Cullinane et al. 1992, offer not only the ability to analyze single hydraulic loading condition, but also capability to analyze looped WDNs, employing the water reliability. This paper, however, brings about a novel and simple, yet very effective computational-method, called Spanning-Tree Matrix Approach (STMA), to model tree-shaped WDNs, rendering the consideration of single hydraulic loading possible.

To the best of authors knowledge, all the studies stated in the literature address only the minimum-nodal-head constraint and, sometimes, the maximum-velocity constraint. Proposed STMA in this study can function, similar to dynamic programming problems, as an adaptable sub-system in an optimization algorithm which utilizes the evaluation of desired hydraulic constraints – minimum-nodal-head, and maximum-pipe-head-loss - for WDNs that have tree-shaped architecture. A case study has also been presented to demonstrate the coupled function of STMA and a modern stochastic optimization algorithm, the HBMO algorithm, that takes into consideration not only the constraint of minimum-nodal-head but also the maximum-pipe-head-loss.

# FORMULATION OF THE OPTIMIZATION MODEL

Optimal WDN design problem focused in this study can be written as,

$$Min\ Z = \sum_{i=1}^{N} C_i(D,L) \tag{1}$$

subjects to the hydraulic design constraints,

$$H_{R\ j} \geq H_{R\ j}^{min} \quad j = 1,2,3, \dots, nd \tag{2}$$

$$g_{FF\ i} \leq g_{FF\ i}^{max} \quad i = 1,2,3, \dots, np \tag{3}$$

where, $Z$= total cost involved with the WDN; $N$ = number of pipes in the water distribution network; $C_i(D,L)$ = cost of the $i^{th}$ pipe having diameter $D$ and length $L$; $H_{R\ j}$ = residual water head available at the $j^{th}$ node; $H_{R\ j}^{min}$ = minimum residual water head required at the $j^{th}$ node; $g_{FF\ i}$ = friction and fitting loss gradient in the $i^{th}$ pipe; $g_{FF\ i}^{max}$ = maximum allowable friction and fitting loss gradient in the $i^{th}$ pipe; $nd$ = number of demand nodes; and $np$ = number of pipes.

The flowchart shown in Fig. 1, depicts the algorithm which accommodates the above optimization model. It is typically required the WDN data, such as reservoir elevation, nodal demands, pipe lengths, nodal elevations, and the pipe connectivity layout. In addition to WDN data, commercially available pipe diameters and their unit costs are needed as the fulfillment of the data requirement.

Once the data are fed to the coupled hydraulic and optimization model, it iteratively seeks for the best solution until the stopping criterion is met. The solution thus obtained can be regarded as the global optimal solution, or rather the optimal combination of pipe diameters that minimizes the cost of the WDN to the maximum extent.

Next, it can be discussed how the proposed network hydraulic simulation by STMA model can be done to analyze a WDN that takes a tree-shaped architecture.

## SPANNING-TREE MATRIX APPROACH

The goal of using Spanning-Tree Matrix Approach (STMA) is to obtain the nodal-residual-head vector $[H_R]$ and the friction and fitting loss gradient vector $[g_{FF}]$ for the selected diameter vector [d] abiding to the algorithm shown in the Fig. 1. To accomplish this task, the basic steps that are followed are listed in the Fig. 2.

In order for the employed mathematical concepts in the computations to be clearly understood, it is convenient to consider a fairly small, tree-shaped, dummy WDN that has six consumer demand nodes. As the flow occurs in one direction only, the spanning-tree of WDN here is a directed graph rooted at the reservoir node (see the Fig. 3).

The reservoir, which is the root of the spanning-tree, is denoted by $N_0$. However, $N_i$s ($i = 1,2,3,4,5,6$) are the consumer demand nodes of the WDN; $L_i$s ($i = 1,2,3,4,5,6$) are the lengths of pipes; $Q_i$s ($i = 1,2,3,4,5,6$) are the flows in pipes in the indicated direction; and $D_j$s ($j = 1,2,3,4,5,6$) are the nodal water demands.

In matrix notation, the nodal water demands and pipeline flows can be written as (4) and (5), respectively.

$$[D] = [D_j] = [D_1 \ D_2 \ D_3 \ D_4 \ D_5 \ D_6]^T$$

$$j = 1,2,3,4,5,6 \tag{4}$$

$$[Q] = [Q_i] = [Q_1 \ Q_2 \ Q_3 \ Q_4 \ Q_5 \ Q_6]^T \tag{5}$$

$$i = 1,2,3,4,5,6$$

Then the flow vector $[Q]$, computed to facilitate the continuity equation, is given by,

$$[Q] = [Q_i] = [T_D^Q] \cdot [D] \tag{6}$$

where, $[T_D^Q]$ is not only a form of adjacency matrix that illustrates the connectivity of pipes through nodes, but also a linear-transformation matrix that transforms nodal demands in to pipeline flows. The matrix $[T_D^Q]$ for the WDN corresponding to the Figure 3 can be considered as (7). Intuitively, for example, if the pipe number 1 has a flow of $Q_1$, then $Q_1$ is equal to the summation of the demands of all the downstream nodes, resulting 1s at the corresponding columns and 0s for all the other entries of the first row of $[T_D^Q]$. Similarly, all the entries of all the other rows of $[T_D^Q]$ can be determined. (7)

$$[T_D^Q] = \begin{bmatrix} 1 & 0 & 1 & 1 & 0 & 0 \\ 0 & 1 & 0 & 0 & 1 & 1 \\ 0 & 0 & 1 & 0 & 0 & 0 \\ 0 & 0 & 0 & 1 & 0 & 0 \\ 0 & 0 & 0 & 0 & 1 & 0 \\ 0 & 0 & 0 & 0 & 0 & 1 \end{bmatrix}$$

The matrix form, or rather the vector form, of the selected solution can be seen as given in (8). There, the $d_i$ is the diameter of the $i^{th}$ pipe.

$$[d] = [d_i] = [d_1 \ d_2 \ d_3 \ d_4 \ d_5 \ d_6]^T \tag{8}$$

$$i = 1,2,3,4,5,6$$

Friction-loss-gradient of the $i^{th}$ pipe, $h_{f\ i}$, is then calculated according to the Hazen-Williams Equation shown in (9).

$$h_{f\ i} = \frac{10.666 * Q_i^{1.85}}{C_{HW}^{1.85} * d_i^{4.87}} \tag{9}$$

$$i = 1,2,3,4,5,6$$

The matrix form of the set of friction-loss gradients for all the pipes can be identified as (10).

$$[h_f] = [h_{f\ i}] = \frac{10.666}{C_{HW}^{1.85}} \cdot [[Q]^{o\ (1.85)} ./ [d]^{o\ (4.87)}] \tag{10}$$

where, $[Q]^{o\ (1.85)}$ is the notation for the vector obtained by the element-wise exponentiation of the elements in $[Q]$ to the power 1.85; and "./" is the operator "Hadamard Division" for denoting the element wise division between two vectors $[Q]^{o\ (1.85)}$ and $[d]^{o\ (4.87)}$.

The equations (6) and (10) can be combined to obtain the equation (11) as shown below.

$$[h_f] = [h_{f\ i}] = \frac{10.666}{C_{HW}^{1.85}} \cdot [([T_D^Q] \cdot [D])^{o\ (1.85)} ./ [d]^{o\ (4.87)}] \tag{11}$$

To account for the fitting losses that occur in each pipe, friction-loss gradient vector is multiplied by the fitting-loss coefficient $C_{ft}$ to obtain the $[g_{FF}]$, the gradient vector for both friction and fitting losses.

$$[g_{FF}] = C_{ft} \cdot [h_f] \tag{12}$$

The pipe-length vector $[L]$ of the WDN shown in Figure 3, can be regarded equal to the expression (13).

$$[L] = [L_i] = [L_1\ L_2\ L_3\ L_4\ L_5\ L_6]^T \tag{13}$$

$$i = 1,2,3,4,5,6$$

The pipeline-head-loss vector $[F]$, is the vector whose $i^{th}$ element is the product of $i^{th}$ elements of $[L]$ and $[g_{FF}]$. The matrix notation of this element-wise multiplication of two vectors, $[L]$ and $[g_{FF}]$, is given by the "Hadamard Product $\odot$" as shown in (14).

$$[F] = [F_i] = [L] \odot [g_{FF}] \tag{14}$$

The equations (12) and (14) can be merged to obtain the matrix equation (15).

$$[F] = [L] \odot C_{ft}.[h_f] \tag{15}$$

As the hydraulic-head-loss along the flow path is a continuous function, total algebraic summation of the head-loss around a closed path, is equal to zero. The mathematical notation of this law of conservation of energy is given by,

$$\sum_{i=1}^{np_L} F_i = 0 \tag{16}$$

where, $npL$ is the number of pipes in a closed path.

As the focus here in this study is directed only towards the tree-shaped WDNs, the adaptation of the above (16) is shown in (17) as a matrix form.

$$[H] = E_0[1] - [T_H].[F] \tag{17}$$

where, $[H]$ is the nodal-head vector; $E_0$ is the elevation of the reservoir node $N_0$ measured from the mean sea level (MSL); $[1]$ is the vector in which each element is equal to 1; $[T_H]$ is the transformation matrix that transforms pipeline-head-loss values $[F]$ into total head-loss values at all the nodes, and total head-loss values are given by $[T_H].[F]$. The matrix expression of $[T_H]$, however, is given by (18).

$$\tag{18}$$

$$[T_H] = \begin{bmatrix} 1 & 0 & 0 & 0 & 0 & 0 \\ 0 & 1 & 0 & 0 & 0 & 0 \\ 1 & 0 & 1 & 0 & 0 & 0 \\ 1 & 0 & 0 & 1 & 0 & 0 \\ 0 & 1 & 0 & 0 & 1 & 0 \\ 0 & 1 & 0 & 0 & 0 & 1 \end{bmatrix}$$

It can be noticed, if $N_4$ is considered for example, head-loss occurs only through pipe number 1 and 4, hence first and fourth column entries of row number four are equal to 1, resulting all the other entries in that row being equal to zero. Likewise, entries of all the other rows of $[T_H]$ can be determined by mere layout of the WDN concerned.

By considering expressions (17) and (15) together, the following expression (19) can be derived for $[H]$.

$$[H] = E_0[1] - [T_H].([L] \odot C_{ft}.[h_f]) \tag{19}$$

The nodal-elevation vector $[E]$ has entries $E_j$s where $E_j$ is the nodal-elevation of the $j^{th}$ node measured from MSL. $[E]$ is given by expression (20).

$$[E] = [E_j] = [E_1 \; E_2 \; E_3 \; E_4 \; E_5 \; E_6]^T \tag{20}$$

$$j = 1,2,3,4,5,6$$

As the residual water head $H_R$ at each node is equal to the difference between the nodal-head and nodal-elevation, the matrix form is given by (21).

$$[H_R] = [H] - [E] \tag{21}$$

By merging the expressions (19) and (21) together, below expression (22) can be obtained.

$$[H_R] = E_0[1] - [T_H].([L] \odot C_{ft}.[h_f]) - [E] \tag{22}$$

As a summary, from the proposed STMA, $[g_{FF}]$ and $[H_R]$ are computed for selected $[d]$. Hence, it is now possible to output constraint evaluations of $[g_{FF}]$ and $[H_R]$ in order to continue on with the algorithm being described under Fig. 1.

## CASE STUDY: WARAPITIYA SERVICE ZONE, SRI LANKA

The general optimization algorithm that includes proposed STMA simulation described under Fig. 1, was implemented for the WDN scheme planned for Warapitiya Service Zone, Sri Lanka, aiming for obtaining the combination of pipe diameters that minimizes the cost function while satisfying the design constraints, as per the optimization model formulation described. In this study the modern stochastic optimization algorithm Senavirathna et al. 2022 have used, i. e., Honey-Bee Mating Optimization Algorithm (HBMO), was updated with the STMA component to facilitate the hydraulic simulations for the tree-shaped WDN. The minimum-allowable-nodal-hydraulic-head value $H_R^{min}$ and the maximum-allowable-gradient for friction-and-fitting losses $g_{FF}^{max}$ were considered to be 10 m and 0.005 m/m, respectively. The Hazen-Williams coefficient, $C_{HW}$, and the fitting loss coefficient, $C_{ft}$, were taken to be 130 and 1.15, respectively, for all the pipes.

Here, in fact, are two main attributes of Warapitiya WDN to be an ideal choice for STMA simulation. On the one hand, WDN of Waraptiya Service Zone has the layout architecture of a spanning-tree, and on the other hand, it has substantial number of demand nodes and pipes that it can be considered as a relatively "large-scale" WDN.

The data requirement, as identified for the optimization algorithm, included unit-costs of commercially available pipe sizes, reservoir water head, nodal-water-demands, pipe lengths, nodal-elevations, and the pipe-connectivity layout. The pipe layout of Warapitiya WDN can be seen from the Fig. 4 as a directed-graph to identify the flow directions. It is to be noted that the symbology used under the Warapitiya Service Zone mean the same as that was utilized in STMA methodology explanation. Pipe length data are given in the Table 1.

Nodal-water demand data and nodal-elevation data, including the datum of reservoir elevation, can be referred to as listed in the Table 2. All these data, excluding unit-costs of commercially available pipes, were collected from the National Water Supply and Drainage Board (NWS&DB), Sri Lanka. Unit-cost data for different pipe sizes were obtained from Mohan and Babu 2010, and shown in the Table 3. Corresponding unit-costs for commercially available pipe sizes used in Warapitiya WDN, however, were estimated by interpolating these data using Lagrange-Interpolation-Method.

According to Lagrange-Interpolation-Method, as the number of data points given in Table 3 is 14 (=13+1), the Lagrange-interpolating-polynomial interpolates 14 data points by a polynomial of degree 13, $f_{13}(x)$, defined as,

$$f_{13}(x) = \sum_{k=0}^{13} L_k(x) f(x_k) \qquad (23)$$

where,

$$L_k(x) = \prod_{\substack{m=0 \\ m \neq k}}^{13} \frac{(x - x_m)}{(x_k - x_m)} \qquad (24)$$

The interpolated unit-costs for commercially available pipe diameters are listed in the Table 4.

## RESULTS

From Table 5, it can be observed the contrast between the optimal design solution obtained by HBMO algorithm and the solution implemented by NWS&DB. Most of the diameter-sizes given by HBMO algorithm were less than or equal to those implemented by NWS&DB, and however, the pipes $P_6$, $P_{10}$, $P_{16}$, $P_{18}$, $P_{19}$, and $P_{20}$ diverge from that observation. Moreover, the total cost involved with the optimal design solution obtained by HBMO algorithm was 100172 units whereas the total cost involved with the solution implemented by NWS&DB was 107588 units, exhibiting lower cost for HBMO solution. The optimal design solution obtained by optimization algorithm highlighted in Fig. 1, therefore, was lower in cost compared to randomly guessed reasonable solutions, for no violation in the design constraints in both cases.

Table 6 and Table 7, respectively, tabulate the results, $[g_{FF}]$ and $[H_R]$, obtained at the STMA simulation stage, corresponding to the global optimal solution obtained by HBMO algorithm. The Table 6 shows that the design constraint (3), $g_{FF\ i} \leq g_{FF\ i}^{max}$ is satisfied for all the pipes, $P_i$s ($i = 1,2,...,24$), in the spanning-tree network by delivering the $g_{FF}$ value of each pipe less than $0.005\ m/m$, whereas the Table 7 shows that the design constraint (2), $H_{R\ j} \geq H_{R\ j}^{min}$ is satisfied for all the nodes, $N_j$s ($j = 1,2,...,24$), in the spanning-tree network by delivering the $H_R$ value of each node greater than $10\ m$. Due to high flow rate in the pipe $P_1$, it is observable that, the $g_{FF}$ has come closer to its upper limit. Nevertheless, the application of the obtained global optimal solution for $[d]$ is a *feasible solution*, as all the design constraints are satisfied for all the pipes and all the nodes in the network.

Although in the Table 6 and the Table 7, it is investigated only the constraints $g_{FF\ i} \leq g_{FF\ i}^{max}$ and $H_{R\ j} \geq H_{R\ j}^{min}$, respectively, all the intermediate attributes of the STMA can be adapted to output the constraint evaluations of desired design criteria to be coupled with the stochastic-heuristic optimization algorithm. For example, it is possible to easily evaluate the potential velocity constraints and even pumping constraints in the STMA hydraulic simulation if the designer intends to consider such constraints for optimization algorithm.

To summarize, the results of the study demonstrate that, HBMO algorithm together with the proposed STMA simulation were able to deliver the combination of pipe diameters that minimizes the cost of the given spanning-tree pipe configuration network while satisfying the design constraints being concerned.

## CONCLUSIONS AND RECOMMENDATIONS

Aiming for deriving a hydraulic model that has adaptable and easy-to-adjust hydraulic governing equations, and can output evaluations of study specific design constraints to stochastic optimization algorithms, this study was carried out to propose and utilize the Spanning-Tree Matrix Approach (STMA) that can model large-scale tree-shaped water distribution networks. Taking the WDN of Warapitiya Service Zone in Sri Lanka as a case study, it was tested the use of the STMA model coupled with the Honey-Bee Mating Optimization (HBMO) algorithm to find the combination of pipe diameters that minimizes the cost of the pipe network. The results show that the STMA is successful in modeling a tree-shaped WDN of given design criteria, and proved its applicability as a hydraulic model in stochastic optimization algorithms. Furthermore, proposed STMA can be updated with desired governing equations or design criteria that are not discussed in this study.

## DATA AVAILABILITY STATEMENT

All data and models used during the study appear in this article itself. However, the codes that support the findings of this study are available from the corresponding author upon reasonable request.

## ACKNOWLEDGEMENTS

Authors are extending sincere thanks to the National Water Supply and Drainage Board, Sri Lanka for their data used in this study.

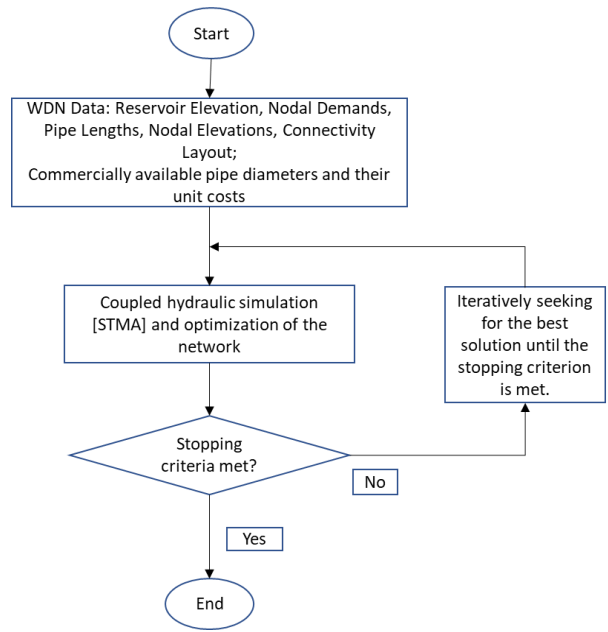

**Fig. 1.** Flowchart of the optimization algorithm considered in this study

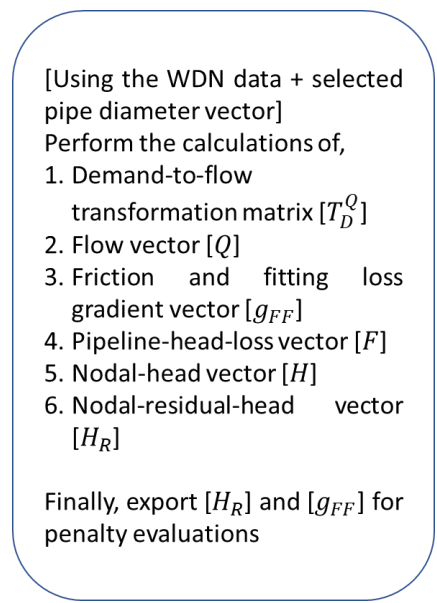

**Fig. 2.** Basic steps involved with the version of the STMA used in the study

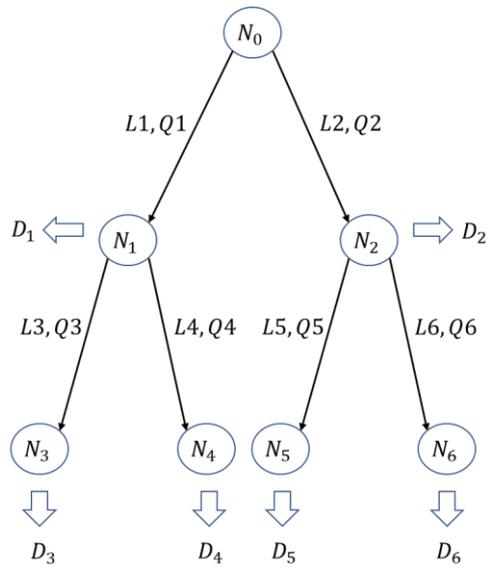

**Fig. 3.** The dummy WDN for the demonstration of the functionality of STMA.

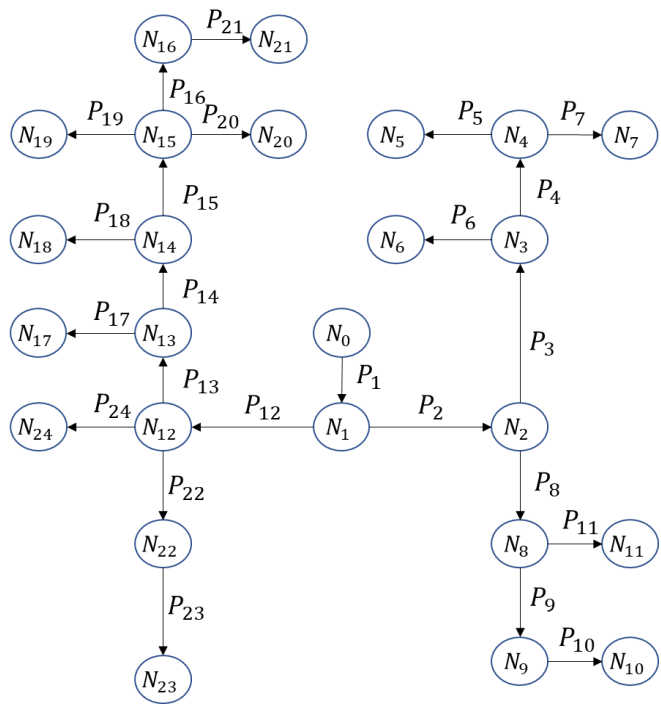

**Fig. 4.** Pipe layout directed-graph for WDN of Warapitiya Service Zone, Sri Lanka.

**Table 1.** Lengths of pipes being deployed for Warapitiya WDN

| Pipe ($P_i$) | Length ($L_i$) / (m) |
|:---:|:---:|
| $P_1$ | 250 |
| $P_2$ | 160 |
| $P_3$ | 740 |
| $P_4$ | 80 |
| $P_5$ | 290 |
| $P_6$ | 350 |
| $P_7$ | 130 |
| $P_8$ | 330 |
| $P_9$ | 160 |
| $P_{10}$ | 180 |
| $P_{11}$ | 410 |
| $P_{12}$ | 30 |
| $P_{13}$ | 260 |
| $P_{14}$ | 630 |
| $P_{15}$ | 300 |
| $P_{16}$ | 530 |
| $P_{17}$ | 850 |
| $P_{18}$ | 520 |
| $P_{19}$ | 290 |
| $P_{20}$ | 620 |
| $P_{21}$ | 140 |
| $P_{22}$ | 300 |
| $P_{23}$ | 1280 |
| $P_{24}$ | 1160 |

**Table 2.** Nodal-water-demands and nodal-elevations of Warapitiya WDN

| Node ($N_j$) | Demand ($D_j$) / ($m^3/day$) | Elevation ($E_j$) / ($m$) |
|---|---|---|
| $N_0$ | - | 506.0 |
| $N_1$ | 3.75 | 485.0 |
| $N_2$ | 381.56 | 475.0 |
| $N_3$ | 292.50 | 455.0 |
| $N_4$ | 148.13 | 460.0 |
| $N_5$ | 48.75 | 455.0 |
| $N_6$ | 11.25 | 450.0 |
| $N_7$ | 91.88 | 460.0 |
| $N_8$ | 101.25 | 480.0 |
| $N_9$ | 45.00 | 480.0 |
| $N_{10}$ | 18.75 | 455.0 |
| $N_{11}$ | 43.13 | 476.0 |
| $N_{12}$ | 71.25 | 482.0 |
| $N_{13}$ | 133.13 | 479.0 |
| $N_{14}$ | 105.00 | 472.0 |
| $N_{15}$ | 116.25 | 462.5 |
| $N_{16}$ | 61.88 | 445.0 |
| $N_{17}$ | 61.88 | 450.0 |
| $N_{18}$ | 28.13 | 474.0 |
| $N_{19}$ | 33.75 | 450.0 |
| $N_{20}$ | 28.13 | 452.0 |
| $N_{21}$ | 20.63 | 450.0 |
| $N_{22}$ | 76.88 | 475.0 |
| $N_{23}$ | 71.25 | 450.0 |
| $N_{24}$ | 58.13 | 450.0 |

**Table 3.** Unit-cost data for corresponding pipe sizes (Mohan & Babu, 2010)

| Diameter / ($mm$) | Unit-cost |
|---|---|
| 25.4 | 2 |
| 50.8 | 5 |
| 76.2 | 8 |
| 101.6 | 11 |
| 152.4 | 16 |
| 203.2 | 23 |
| 254.0 | 32 |
| 304.8 | 50 |
| 355.6 | 60 |
| 406.4 | 90 |
| 457.2 | 130 |
| 508.0 | 170 |
| 558.8 | 300 |
| 609.6 | 550 |

**Table 4.** Commercially available pipe sizes and their corresponding interpolated unit-costs

| Diameter / ($mm$) | Interpolated unit-cost |
|---|---|
| 55 | 5.0259 |
| 79 | 8.4781 |
| 97 | 10.6801 |
| 140 | 14.0679 |
| 198 | 22.5046 |
| 246 | 29.6739 |

**Table 5.** Comparison between the solution obtained by HBMO+STMA and the guessed solution implemented by NWS&DB.

| Pipe $P_i$ | Optimal Solution [d] by HBMO+STMA / $(mm)$ | Guessed solution implemented by NWS&DB / $(mm)$ |
|---|---|---|
| $P_1$ | 198 | 246 |
| $P_2$ | 198 | 198 |
| $P_3$ | 198 | 198 |
| $P_4$ | 97 | 140 |
| $P_5$ | 79 | 79 |
| $P_6$ | 79 | 55 |
| $P_7$ | 79 | 79 |
| $P_8$ | 97 | 140 |
| $P_9$ | 79 | 79 |
| $P_{10}$ | 97 | 55 |
| $P_{11}$ | 55 | 79 |
| $P_{12}$ | 198 | 198 |
| $P_{13}$ | 140 | 198 |
| $P_{14}$ | 140 | 140 |
| $P_{15}$ | 97 | 140 |
| $P_{16}$ | 97 | 79 |
| $P_{17}$ | 79 | 79 |
| $P_{18}$ | 97 | 55 |
| $P_{19}$ | 79 | 55 |
| $P_{20}$ | 79 | 55 |
| $P_{21}$ | 55 | 55 |
| $P_{22}$ | 79 | 97 |
| $P_{23}$ | 55 | 79 |
| $P_{24}$ | 55 | 79 |

**Table 6.** The constraint (3), pipe-constraint [$g_{FF\ i} \leq g_{FF\ i}^{max}$] satisfaction demonstration for the solution obtained by HBMO

| Pipe $P_i$ | Optimal Solution [d] by HBMO / $(mm)$ | $g_{FF\ i}$ / $(m/m)$ | Design Constraint $g_{FF\ i} \leq g_{FF\ i}^{max}$ satisfied? |
|---|---|---|---|
| $P_1$ | 198 | 0.0040 | Yes |
| $P_2$ | 198 | 0.0014 | Yes |
| $P_3$ | 198 | 0.0004 | Yes |
| $P_4$ | 97 | 0.0034 | Yes |
| $P_5$ | 79 | 0.0003 | Yes |
| $P_6$ | 79 | 0.0000 | Yes |
| $P_7$ | 79 | 0.0011 | Yes |
| $P_8$ | 97 | 0.0019 | Yes |
| $P_9$ | 79 | 0.0006 | Yes |
| $P_{10}$ | 97 | 0.0000 | Yes |
| $P_{11}$ | 55 | 0.0016 | Yes |
| $P_{12}$ | 198 | 0.0008 | Yes |
| $P_{13}$ | 140 | 0.0021 | Yes |
| $P_{14}$ | 140 | 0.0010 | Yes |
| $P_{15}$ | 97 | 0.0028 | Yes |
| $P_{16}$ | 97 | 0.0003 | Yes |
| $P_{17}$ | 79 | 0.0005 | Yes |
| $P_{18}$ | 97 | 0.0000 | Yes |
| $P_{19}$ | 79 | 0.0002 | Yes |
| $P_{20}$ | 79 | 0.0001 | Yes |
| $P_{21}$ | 55 | 0.0004 | Yes |
| $P_{22}$ | 79 | 0.0027 | Yes |
| $P_{23}$ | 55 | 0.0004 | Yes |
| $P_{24}$ | 55 | 0.0028 | Yes |

**Table 7**. The constraint (2), node-constraint $[H_{R\ j} \geq H_{R\ j}^{min}]$ satisfaction demonstration for the solution obtained by HBMO

| Node $N_j$ | $H_{R\ j}$ / (m) | Design Constraint $H_{R\ j} \geq H_{R\ j}^{min}$ satisfied? |
|---|---|---|
| $N_1$ | 20.0088 | Yes |
| $N_2$ | 29.7801 | Yes |
| $N_3$ | 49.4855 | Yes |
| $N_4$ | 44.2133 | Yes |
| $N_5$ | 49.1135 | Yes |
| $N_6$ | 54.4775 | Yes |
| $N_7$ | 44.0688 | Yes |
| $N_8$ | 24.1675 | Yes |
| $N_9$ | 24.0770 | Yes |
| $N_{10}$ | 49.0731 | Yes |
| $N_{11}$ | 27.5114 | Yes |
| $N_{12}$ | 22.9846 | Yes |
| $N_{13}$ | 25.4312 | Yes |
| $N_{14}$ | 31.7940 | Yes |
| $N_{15}$ | 40.4496 | Yes |
| $N_{16}$ | 57.7719 | Yes |
| $N_{17}$ | 53.9765 | Yes |
| $N_{18}$ | 29.7702 | Yes |
| $N_{19}$ | 52.8990 | Yes |
| $N_{20}$ | 50.8724 | Yes |
| $N_{21}$ | 52.7147 | Yes |
| $N_{22}$ | 29.1779 | Yes |
| $N_{23}$ | 48.9941 | Yes |
| $N_{24}$ | 51.7607 | Yes |